# The Generalized Word Problem for Braid Groups

Elie Feder
Kingsborough Community College
efeder@kbcc.cuny.edu

**Abstract**: The braid group has recently attracted much attention. This is primarily based upon the discovery of its usage in various cryptosystems [AAG],[KLCHKP]. One major focus of current research has been in solving decision problems in braid groups and analyzing their complexity. Among the decision problems which have been investigated are: the word problem, the conjugacy problem and the conjugacy search problem. One decision problem which has not been thoroughly pursued is the generalized word problem. It can be shown that this problem is, in general, unsolvable for braid groups of index greater then or equal to 5. However, we will show that for the special case of cyclic subgroups, this problem is, in almost all cases, solvable. We will design an algorithm, implement it on a multitape Turing machine, and analyze its complexity.

## 1 Introductory Remarks

Throughout this paper, we will deal with two different presentations of the braid group: the Artin presentation and the BKL presentation.
Artin defined $B_n$, the braid group with braid index n, by the following generators and relations:
Generators: $\sigma_1$, $\sigma_2$, …, $\sigma_{n-1}$
Defining relations: (1) $\sigma_s \sigma_t = \sigma_t \sigma_s$, for $|t-s|>1$
(2) $\sigma_s \sigma_t \sigma_s = \sigma_t \sigma_s \sigma_t$, for $|t-s|=1$.

An alternate presentation for $B_n$ was suggested by Birman, Ko and Lee [BKL]. This presentation was introduced in order to reduce the complexity of decision problems in the braid group. It consists of the following generators and relations:
Generators: $a_{ts}$, $1 \leq s < t \leq n$
Defining relations: 1) $a_{ts} a_{rq} = a_{rq} a_{ts}$  if $(t-r)(t-q)(s-r)(s-q)>0$
2) $a_{ts} a_{sr} = a_{tr} a_{ts} = a_{sr} a_{tr}$
Notice that $\sigma_i = a_{(i+1,i)}$. Thus, the Artin generators are a subset of the BKL generators.

In studying the generalized word problem, we will refer to the solution of the word problem and its complexity. Thus, let us say a word about the word problem in the braid group. This problem could be formulated as follows: Given a word W in $B_n$, can we decide in a finite number of steps whether or not W is equivalent to the identity element in $B_n$. This is equivalent to determining if two words $W_1$ and $W_2$ are equivalent. For, to say that $W_1 = W_2$, is equivalent to saying that $W_1 W_2^{-1} = 1$.
There have been a number of solutions for the word problem in the braid group. The main idea of several of these approaches is to find an algorithm for putting a braid into a unique normal form and thereby compare any two braids. This was first introduced by Garside[G] in 1969. His solution is exponential in braid index and word length. In 1992, Thurston[T] improved upon Garside's approach and proved that for word length $|W|$, there exists an algorithmic solution to the word problem in $B_n$ with complexity

$O(|W|^2 n \log n)$. In 2001, [BKL] further improved upon the algorithm of Garside-Thurston by introducing the BKL-presentation for the braid group. Using their new generators and defining relations, they mirrored the algorithm of Garside-Thurston and reduced the complexity of the word problem in the braid group to $O(|W|^2 n)$. Although there have been a number of other approaches to the solution to the word problem in the braid group, for our purposes these will suffice.

## 2 The Generalized Word Problem

We will deal with the generalized word problem for braid groups. The problem is formulated as follows: Given a braid group $B_n$ and given H, a subgroup of $B_n$, can we decide in a finite number of steps whether or not a given word W in $B_n$ is also contained in H? Notice that the word problem is the special case of the generalized word problem when H is equal to the trivial subgroup in $B_n$.

A general solution to this problem is not possible. The reasoning is as follows. In [BD], it is shown that $B_5$ contains an isomorphic copy of the direct product of two free groups of rank 2, as a subgroup. Additionally, in [Mih], it was shown that the direct product of two free groups has unsolvable generalized word problem. This, being the case, one can conclude that the generalized word problem for braid groups of index greater then or equal to 5 is unsolvable. However, we will prove that for the specific case of cyclic subgroups whose generator has nonzero exponential sum, the generalized word problem is, in fact, solvable. We will design an algorithm and implement it on a multitape Turing machine, and thereby analyze the complexity of this problem.

Being that we will only consider the generalized word problem for the case of cyclic subgroups, we will rephrase the problem in the following equivalent form: Given a braid x and a braid y, both in $B_n$, can we find an algorithm to determine whether y is in the cyclic subgroup generated by x or not? If it is, can we find the exponent k such that $x^k = y$, and demonstrate this equivalence? What is the running time for this algorithm? One can think of this problem as a kind of logarithm problem in the context of braid groups. Namely, evaluate k where $k = \log_x y$.

Note: This problem was raised in the newsgroup sci.crypt[S] and attracted much interest.

## 3 Some Preliminary Lemmas

To solve this problem, we introduce a definition and some lemmas.

<u>Definition</u>: We define a function $\exp: B_n \to \mathbf{Z}$ where for any $\beta \in B_n$, $\exp(\beta)$=the sum of the exponents of $\beta$.
<u>Example</u>: $\exp(\sigma_1 \sigma_3^{-3} \sigma_2^2 \sigma_1) = 1-3+2+1 = 1$.

<u>Lemma 1</u>: exp is invariant under braid equivalence. In other words if u=v in $B_n$, then $\exp(u)=\exp(v)$.
<u>Proof</u>: Let u=v in $B_n$. Then we can get from u to v by a chain of elements in $B_n$, where each element in the chain can be attained from the previous element by applying one of the defining relations in $B_n$. Thus, it suffices to show that applying any of the defining relations in $B_n$ will not change the value of exp. Let us, therefore, consider the defining relations in $B_n$.
1. $\sigma_i\sigma_j=\sigma_j\sigma_i$ where $|i-j|>1$. Well, $\exp(\sigma_i\sigma_j) = 2 = \exp(\sigma_j\sigma_i)$.
2. $\sigma_i\sigma_{i+1}\sigma_i=\sigma_{i+1}\sigma_i\sigma_{i+1}$ .Well, $\exp(\sigma_i\sigma_{i+1}\sigma_i) = 3 = \exp(\sigma_{i+1}\sigma_i\sigma_{i+1})$.

Note: If we have an occurrence of $\sigma_i\sigma_i^{-1}$ or $\sigma_i^{-1}\sigma_i$ and we replace it with 1, we also do not alter the value of exp. This is because $\exp(\sigma_i\sigma_i^{-1}) = \exp(1) = \exp(\sigma_i^{-1}\sigma_i) = 0$. The lemma is thus proven.

Lemma 2: exp is additive. Namely $\exp(u\cdot v)=\exp(u)+ \exp(v)$
Proof: Since $u\cdot v$ is simply the concatenation of the words u and v, the sum of the exponents of the new word $u\cdot v$ is simply the sum of the exponents of u plus the sum of the exponents of v. Any cancellation caused by the concatenation is simply the removal of $\sigma_i\sigma_i^{-1}$ or $\sigma_i^{-1}\sigma_i$, both of whose exponent sum is zero. Thus removal of these terms will not affect the total exponent sum.

Lemma 3: $\exp(\beta^{-1})= -\exp(\beta)$
Proof: Let $\beta= \sigma_{i1}^{p_1}\sigma_{i2}^{p_2}\ldots\sigma_{in}^{p_n}$. Then $\exp(\beta)= p_1+p_2+\ldots+p_n$. Now $\beta^{-1}= \sigma_{in}^{-p_n}\ldots\sigma_{i2}^{-p_2}\sigma_{i1}^{-p_1}$. Then $\exp(\beta^{-1})= -p_1-p_2-\ldots-p_n= -(p_1+p_2+\ldots+p_n)= -\exp(\beta)$.

Lemma 4: If $y=x^k$, then $\exp(y)= k\cdot\exp(x)$.
Proof: This follows directly from the additivity of exp (lemma 2).

Notice that all of these lemmas hold true if we represent a braid in the BKL-generators instead of in the Artin generators. This is the case because all of the relations in the BKL-presentation preserve the value of exp, just as the relations in the Artin presentation.

## 4 The Solution

a) If we are given that $y=x^k$ for some k, and if we further assume that $\exp(x)\neq 0$, then we can easily solve for k. Notice that $\exp(y)= k\cdot\exp(x)$, by lemma 4. Thus, $k= \exp(y)/\exp(x)$.
b) Now, let us assume that we are given x, y in $B_n$, with $\exp(x)\neq 0$ and we are told to determine if $y=x^k$ for some k. In other words, is y in the cyclic subgroup generated by x? We proceed as follows. Assume it is. Then $y=x^k$, for some k. We can compute the only possible value for k. Namely, $k= \exp(y)/ \exp(x)$. Now compute $x^k$. Decide if $y=x^k$ by the known algorithms for the word problem[CKLHC].
Note: The above solution is limited to the case where $\exp(x)\neq 0$. For any braid word x, where $\exp(x)=0$, a different solution would be required. As $\exp(x)=0$ for any $x\in [B_n,B_n]$, the commutator of $B_n$, such a solution would certainly be desirable.

## 5 The Outline of the Algorithm

Given: $x,y \in B_n$.
Question: Is y in the cyclic subgroup generated by x?

Step 1: Compute exp(x) and exp(y)
Step 2: Compute exp(y)/exp(x)=k
Step 3: Compute $x^k$
Step 4: Decide $x^k$=y? If yes, then y is in the cyclic subgroup generated by x and we know the power, k. If no, then y is not in the cyclic subgroup generated by x.
<u>Note</u>: In the case that exp(y)=0, then k=0. By convention, define $x^0$=1. Thus, this reduces to the word problem, as a subcase.

## 6  Programming Language for Multitape Turing Machine

We will follow the model of papers by Domanski and Anshel [Do1,Do2,AD] which use a programming language for multitape-Turing machines. We assume a familiarity with the basic model of a multitape-Turing machine which consists of an input tape, k worktapes, and an output tape. Let the symbol being read by tape head A on tape X be denoted as X(A). Moving tape head A on X one tape cell to the right will be denoted by A=A+1. Similarly, moving tape head A on X one tape cell to the left will be denoted by  A=A-1. We allow standard high-level programming constructions such as IF-THEN-ELSE, DO-WHILE and DO-UNTIL… One can verify that these constructions can be converted to standard Turing machine notation. We also allow the number theoretic operations mod and div.  In the course of the algorithm, we will denote the beginning of an explanatory comment with  /* ,  and the end of such a comment with  */.

## 7  The Formal Algorithm GWP(X,Y)

<u>Input</u>: X,Y$\in B_n$ , exp(X)$\neq$0.
<u>Task</u>: Determine if there exists an integer c such that $X^c$=Y.
<u>Output</u>: GWP(X,Y)=  "$X^c$=Y", if  Y is the c-th power of X
"Y is not a power of X", if there does not exist an integer c such that $X^c$=Y.

We will use a multitape-Turing machine with three worktapes; X,Y and Z.
Tape X will contain the braid word X=$X_1 X_2 \ldots X_n$  followed by the $ symbol.
Tape Y will contain the braid word Y=$Y_1 Y_2 \ldots Y_m$  followed by the $ symbol.
$X_1,\ldots,X_n,Y_1,\ldots,Y_m$ each denote a generator of the braid group or its inverse.
We will define sign($X_i$)=1 if $X_i$ is a generator and sign($X_i$)=-1 if $X_i$ is the inverse of a generator.
Tape Z will contain the integers.
A will be a tape head situated on tape X and begins under the first letter of X, namely, $X_1$.
B will be a tape head situated on tape Y and begins under the first letter of Y, namely $Y_1$.
EX and EY will be tape heads situated on tape Z and will begin under zero.
/* EX and EY will be used to keep track of the exponent sum of X and Y */

Do Until X(A)=$
    If sign(X(A))=1,
        Then Do;
            EX=EX+1

```
                End;
        Else Do;
                EX=EX-1
            End;
        A=A+1
End; /* of DO-UNTIL */
/* Z(EX) is now equal to exp(X), the exponent sum of X. */
Do Until Y(B)=$
        If sign(Y(B))=1,
                Then Do;
                        EY=EY+1
                End;
        Else Do;
                EY=EY-1
            End;
        B=B+1
End  /* of DO-UNTIL */
/* Z(EY) is now equal to exp(Y), the exponent sum of Y. */

If Z(EY) mod Z(EX) =0,  /* exp(X) divides exp(Y) */
        Then Do
                c= Z(EY) div (Z(EX))
                Form U=X^c
                /* U is formed by concatenating c copies of X. This can be done in linear
                time on a multitape Turing machine as is done in [Do2] */
                /* Comparison(G,H) is assumed to be a routine which, given G,H∈ $B_n$,
                determines if G=H under the relations of $B_n$. The output is True if G=H
                and False if G≠H   See [CKLHC]. If we assume that L is the minimal of
                the canonical lengths of G and H, then this routine runs in $O(L^2 n \log n)$ if
                G and H were given in Artin generators and in $O(L^2 n)$ if G and H were
                given in BKL generators. */
                If Comparison(U,Y)=True
                        Then Do
                                Output "$X^c$=Y".
                        End;
                Else Do;
                        Output "X is not a power of Y"
                End;
Else Do;
        Output "X is not a power of Y"
    End;
End;
```

## 8  Complexity Analysis of GWP(X,Y)

We will let M be the maximal length of X and Y. The computation of both Z(EX) and Z(EY) is bounded by M. Forming $X^c$ is linear in the length of $X^c$. What is the length of $X^c$? Well, c is bounded by M and, therefore, the length of $X^c$ is certainly bounded by $M^2$. As mentioned above(solution to word problem), if X and Y are given in BKL generators, then the routine Comparison(G,H) is bound by $O(L^2n)$, where L is the minimal length of G and H( $O(L^2 n \log n)$ if we use the Artin generators). In our case where we are comparing $X^c$ and Y, this minimal length is bound by M because the length of Y is bound by M. Thus, Comparison($X^c$,Y) is bound by $O(M^2n)$. Thus, the entire algorithm GWP(X,Y) is bound by $O(M^2n)$ ( $O(M^2 n \log n)$ if we use the Artin generators). It is interesting to compare the complexity of the word problem, determining if X=Y, and the instance of the generalized word problem, determining if Y=$X^c$, for some c. The word problem is bound by $O(L^2n)$, where L is the *minimal* length of X and Y, while the generalized word problem is bound by $O(M^2n)$, where M is the *maximal* length of X and Y.

## 9 Conclusion

We can now conclude with the following theorem:
<u>Theorem 1</u>: The generalized word problem in the case of a cyclic subgroup of the braid group $B_n$ is solvable in almost all cases. More specifically, if we are given
X,Y $\in B_n$, with exp(X)≠0, and if M is the maximal length of X and Y in BKL generators, then we can determine in $O(M^2n)$ whether Y is in the cyclic subgroup generated by X. We can also determine the specific power c such that $X^c$=Y.

Notice that the only properties of braid groups that we used are: 1) that the braid group is exp-invariant( the defining relations do not change the exponent sum) and 2) that the word problem for the braid group is solvable. We can thus formulate the following more general theorem:
<u>Theorem 2:</u> Let G be a group with an exp-invariant presentation and with solvable word problem. Let x be any element in G with nonzero exponent sum. Then the generalized word problem for the case of the cyclic subgroup generated by x is solvable.

## 10 Further Research

1) A problem which is very closely related to this paper is the root problem in the braid group. In a sense, this problem is the opposite of ours. The root problem is as follows: Given a braid word w$\in B_n$, and a positive integer r, we are asked to determine if w has any $r^{th}$ roots. In other words, does there exist a braid b$\in B_n$, such that $b^r$=w? There have been a few results regarding this problem. In [Go] it was shown that any $r^{th}$ roots of a braid w are conjugate to each other. In [St], it was shown that the root problem is solvable. However, no complexity bound was found for this solution. In [Si], this solution was generalized to small Gaussian groups, a generalization of braid groups which was introduced in [DP]. It would be nice if we could find a bound on the complexity of these algorithms.

2) In this paper we have found a quick solution to the generalized word problem in the case of cyclic subgroups of the braid group where exp of the generator is nonzero. Can we find an analogous solution for the case of cyclic subgroups of the braid group which are generated by an element whose exp is zero?

3) A related problem is the word problem for the cyclic amalgamation of braid groups. This problem makes use of the solution to the generalized word problem for the case of cyclic subgroups of the braid group. We, thus, hope to design an algorithm for this word problem and analyze its complexity, based upon the results of this paper.

## Acknowledgements

This paper is an adaptation of a part of my thesis which was written with the help of my advisor Michael Anshel. I thank him for all his help.